\documentclass{elsart3-1}

\input{ot2enc.def}

\usepackage{amssymb}

\def\cC{{\mathcal C}}

\def\cO{{\mathcal O}}

\def\ord{{\rm ord}}

\def\N{{\bf N}}
\def\R{{\bf R}}

\usepackage[francais,russian,english]{babel}

\newtheorem{theorem}{Theorem}[section]

\newtheorem{e-proposition}[theorem]{Proposition}

\newtheorem{e-definition}[theorem]{Definition\rm}

\newtheorem{theoreme}{Th\'eor\`eme}[section]
\newtheorem{lemme}[theoreme]{Lemme}

\newtheorem{definition}[theoreme]{D\'efinition\rm}

\renewcommand{\theequation}{\arabic{equation}}
\def\og{\leavevmode\raise.3ex\hbox{$\scriptscriptstyle\langle\!\langle$~}}
\def\fg{\leavevmode\raise.3ex\hbox{~$\!\scriptscriptstyle\,\rangle\!\rangle$}}

\begin{document}

\begin{frontmatter}
\selectlanguage{francais} \title{La borne de Jacobi pour une diffi{\'e}t{\'e}
d{\'e}finie par un syst{\`e}me quasi r{\'e}gulier}

\vspace{-2.6cm}
\selectlanguage{english}
\title{Jacobi's bound for a diffiety defined by a quasi-regular system}



\author[FO]{F. Ollivier}
\ead{francois.ollivier@lix.polytechnique.fr}
\author[BS]{B. Sadik}
\ead{sadik@ucam.ac.ma}

\address[FO]{ALIEN, INRIA Futurs \&\ LIX, UMR CNRS 7161, {\'E}cole
polytechnique, 91128 Palaiseau CEDEX}
\address[BS]{D{\'e}partement de Math{\'e}matiques,
Facult{\'e} des Sciences Semlalia, B.P. 2390, Avenue Safi, Marrakech,
Maroc}


\begin{abstract}We show that Jacobi's bound for the order of a system
  of ordinary differential equations stands in the case of a diffiety
  defined by a quasi-regular system. We extend the result when there are less
equations than variables and characterize the case when the bound is
reached.

\vskip 0.5\baselineskip

\selectlanguage{francais}
\noindent{\bf R\'esum\'e} On montre que la borne de Jacobi pour
l'ordre d'un syst{\`e}me d'{\'e}quations diff{\'e}rentielles ordinaires est vraie
dans le cas d'une diffi{\'e}t{\'e} d{\'e}finie par un syst{\`e}me quasi r{\'e}gulier.
Nous {\'e}tendons le r{\'e}sulat au
cas o{\`u} il y a moins d'{\'e}quations que d'inconnues et montrons que la
non-nullit{\'e} du jacobien tronqu{\'e} est une condition n{\'e}cessaire et
suffisante pour que la borne soit atteinte.
\vskip 0.5\baselineskip
\noindent

\end{abstract}
\end{frontmatter}

\selectlanguage{english}
\section*{Abridged English version}
In \cite{Jacobi1,Jacobi2}, Jacobi has introduced a bound on the order
of a system of $m$ ordinary differential equations in $m$
unknowns. Let $a_{i,j}$ be the order of the $i^{\rm th}$ equation in
the $j^{\rm th}$ unknown function and $J=\max_{\sigma\in S_{m}}
\sum_{i=1}^{m} a_{i,\sigma(i)}$. He claims that the order of the system
is bounded by $J$. The bound is still conjectural in the general
case. In the setting of differential algebra, Ritt in \cite{Ritt35b}
proved it for linear systems, Lando in \cite{Lando1970} for order 1
systems and Kondratieva {\it et al.} in \cite{Kondratieva1982b} under
the quasi-regularity hypothesis.

Jacobi gave an algorithm to compute the bound in polynomial time,
instead of trying the $m!$ permutations. It has been forgotten and
rediscovered by Kuhn in 1955 (\cite{Kuhn55}), using Egerv{\'a}ry's
results (see \cite{Schrijver05} for historical details).  The idea is
to find a {\em canon}, i.e. integers $\lambda_{1}, \ldots,
\lambda_{m}$ such that, in the matrix $(a_{i,j}+\lambda_{i})$ on can
select maximal entries in each column that are located in different
rows. Jacobi's algorithm computes the unique canon with minimal
$\lambda_{i}$. In the case of $r<m$ equations, we need to compute
$J=\max_{\sigma\in S_{r,m}} \sum_{i=1}^{r} a_{i,\sigma(i)}$, where
$S_{r,m}$ denotes the set of injections $[1,r]\mapsto[1,m]$. That may
be done by completing the matrix $(a_{i,j})$ by $m-r$ rows of $0$ in
order to make it square. Those rows correspond to the orders of $m-r$
generic linear equations of order $0$.

We prove the bound in the context of diffiety extensions (see
\cite{vinogradov,zharinov}) under the same quasi-regularity hypothesis
as in \cite{Kondratieva1982b}. A {\em diffiety} is a real variety of
denumerable dimension equipped with a derivation. {\em Diffiety
morphisms} are $\cC^{\infty}$ mappings commuting with the
derivation. A {\em diffiety extension} $V/U$ is a couple of diffiety
with a canonical surjective projection $\pi:V\mapsto U$ that is a
diffiety morphism. {\em Extension morphisms} also commute with
projections. The {\em trivial extension} of differential dimension
$n$, denoted by $T^{n}/U_{\delta_{0}}$, is defined in coordinates by a
derivation $\delta:=\delta_{0}+\delta_{1}$ with
$\delta_{1}=\sum_{i=1}^{n}\sum_{k=0}^{\infty}
x_{i}^{(k+1)}\partial/\partial x_{i}^{(k)}$. The {\em differential
dimension} of a diffiety extension $V/U$ is the maximal $n$ such that
there is an extension morphism from $V/U$ onto an open set of
$T^{n}/U$. If an extension is of differential dimension $0$, its order
is the dimension of $\pi^{-1}(u\in U)$. For positive dimension, {\em
the order in the coordinates $x_{i}$ at point $(v,u)$} is the maximal
dimension in the neighbourood of $(v,u)$ of the intersection of $V/U$
by differential hyperplanes defined by order zero equations
$h(x)=0$. In may depend of the chosen point. A set of equations
$x_{i}^{(\gamma_{i})}=f_{i}(x)$ defining a diffiety, where the $f_{i}$
depend only of derivatives smaler than $x_{i}^{(\gamma_{i})}$ for some
admissible ordering is said to be a {\em normal form} of the diffiety
in the coordinates $x$.

\begin{e-definition} A system $g(x,y)=0$ of $r$ equations in $m$
unknowns is {\em quasi-regular} at some point $(v,u)\in T^{m}/U$ if for all
$s\in\N$ the jacobian matrix for the derivatives $g_{i}^{(\ell)}$,
$0\le \ell\le s$ with respect to all the derivatives of the unknowns
$x$ appearing in them is of maximal rank $r(s+1)$.
\end{e-definition}

\begin{e-definition} \label{order} Let $g(x,y)=0$ be a system of $r$ equations in $m$
unknowns defining $V$ on some neighbourhood of a point $(v,u)$. We
denote by $\ord_{V,(v,u),x_{j}} g_{i}$ the greatest integer $k$ such
that $\partial g_{i}/\partial x_{j}^{(k)}\neq0$ on every neighbourhood
of $(v,u)$ in $V$, or $-\infty$ if $\partial g_{i}/\partial
x_{j}^{(k)}=0$, $\forall k\in\N$ on some neighbourhood of $(v,u)$ in
$V$. Let $\lambda_{i}$ be a canon for the order matrix $(a_{i,j)}$,
with $a_{i,j}=\ord_{V,(v,u),x_{j}} g_{i}$, $\Lambda=\max_{i}
\lambda_{i}$, $\alpha_{i}=\Lambda -\lambda_{i}$ and
$\beta_{j}=\max_{i} a_{i,j}-\alpha_{i}$. The {\em truncated jacobian
matrix} $\nabla_{V,(v,u)}$ is the matrix $\left(\partial
g_{i}/\partial x_{j}^{(\alpha_{i}+\beta_{j})}\right)$.
\end{e-definition}

\begin{theorem}
Let $g(x,y)$ be a system defining a diffiety extension $V/U$ as a
subdiffiety of the trivial extension $T^{m}/U$ on some neighbourhood
of $(v,u)\in V$ and $J_{V,(v,u)}=\max_{\sigma\in S_{r,m}}
\sum_{i=1}^{r} a_{i,\sigma(i)}$.

i) If $g$ is quasi-regular at $(v,u)\in V$, $V/U$ is of differential
dimension $m-r$ and there exists an open set $\cO$ whose adherence is a
neighbourood of $(v,u)$ and whose every point admits a neighbourood of
order at most $J_{V,(v,u)}$ and possessing a normal form in coordinates $x_{i}$.

ii) If $|\nabla_{V,(v,u)}|\neq0$ at $(v,u)\in V$, then $g$ is
quasi-regular, the differential dimension of $V/U$ is $m-r$, and there is a
neighbourhood having a normal form and of order $J_{V,(v,u)}$ in coordinates
$x_{i}$.

iii) If $g$ is quasi-regular at $(x,u)$ and
$|\nabla_{V,(x,u)}|=0$ on a neighbourhood of $(x,u)$ in $V$, then
there exists an open set $\cO$ whose adherence is a neighbourhood of
$(v,u)$ and whose every point admits a neighbourhood possessing a
normal form and of order strictly lower than $J_{V,(v,u)}$ in coordinates
$x_{i}$.
\end{theorem}

\textsc{Sketch of the proof.} --- The proof follows the scheme
described by Jacobi in \cite{Jacobi2} for the computation of a normal
form. First, we notice that
$\ord_{V,(u,v),x_{j}}g_{i}^{(s)}=\ord_{V,(u,v),x_{j}}g_{i}+s$, using the
commutation rule $[\partial/\partial
x_{j}^{(k+1)},\delta]=\partial/\partial x_{j}^{(k)}$.

(ii) Assume that $|\nabla_{V,(v,u)}|\neq0$ at $(v,u)\in V$. We may
reorder the equations $g_{i}$ by increasing $\alpha_{i}$ and the
$x_{i}$ so that the principal minors made of the first $i$ rows and
columns of $\nabla_{V,(v,u)}$ have a non vanishing determinant $D_{i}$
for all $1\le i\le r$. The jacobian matrix $J_{s}$ of the system
$G_{s}:=\{g_{i}^{(k-\alpha_{i})}|0\le
k\le s\}$ contains a square
submatrix of maximal size obtained by derivating each $g_{i}^{(k-\alpha_{i})}$
with respect to its principal derivatives
$x_{j}^{(\alpha_{i}+k)}$ such that
$\alpha_{j}\le k$. It is block-wise triangular and its
determinant is equal to $\prod_{k=0}^{s} D_{\max\{i| \alpha_{i}\le
k\}}$. It is thus of maximal rank, which proves quasi-regularity. Using
the implicit function theorem for $s$ great enough, we get, using the
fact that $g$ defines $V$ on some neighbourhood of $(v,u)$, a normal form
$x_{j}^{(\alpha_{i}+\beta_{j})}=f_{i}(v)$, $1\le i\le r$. The order in
coordinates $x_{i}$ is then $\sum_{i=1}^{r}
\alpha_{i}+\beta_{i}=J$.

We prove (i) and (iii) by induction on $J_{V,(v,u)}$. In the case $J=0$,
quasi-regularity implies that $\nabla_{V,(v,u)}$ is of full rank, so
the result stands by (ii). If $\nabla_{V,(v,u)}$ is not of full rank
in some neigborough of $(v,u)$, then we may build an equivalent system
with a strictly smaller $J_{V,(v,u)}$ and use (i) recursively. If not, we
consider the open set $\cO'$ where $\nabla_{V,(v,u)}$ is of full rank,
and the interior $\cO''$ of the closed set where $\nabla_{V,(v,u)}$ is
not of full rank. Using (ii) on $\cO'$ and (iii) on $\cO''$, we get (i)
on the open set $\cO'\cup\cO''$.
 
\selectlanguage{francais}
\section{Introduction}
\label{intro}

Dans deux articles posthumes \cite{Jacobi1,Jacobi2}\footnote{On
trouvera la traduction de ces textes sur la page web {\tt
http\char'072//www.lix.polytechnique.fr/\char'176ollivier/JACOBI/jacobi.htm}},
Jacobi a expos{\'e} une borne sur l'ordre d'un syst{\`e}me
d'{\'e}quations diff{\'e}rentielles ordinaires. Celle-ci demeure
conjecturale dans le cas g{\'e}n{\'e}ral. Dans le cadre de
l'alg{\`e}bre diff{\'e}rentielle, elle a {\'e}t{\'e} prouv{\'e}e par
Ritt \cite{Ritt35b} dans le cas lin{\'e}aire, par Lando
\cite{Lando1970} dans le cas d'un syst{\`e}me d'ordre $1$ et par
Kondratieva {\it et al.} sous l'hypoth{\`e}se naturelle de
quasi-r{\'e}gularit{\'e}. Signalons aussi la preuve de cette borne
pour des syst{\`e}mes aux diff{\'e}rences r{\'e}cemment fournie par
Hrushovski dans \cite{Hrushovski2004}.

Jacobi a {\'e}galement d{\'e}crit un algorithme permettant de calculer en
temps polynomial cette borne sans devoir rechercher le maximum de $n!$
combinaisons possibles. Cette contribution a {\'e}t{\'e} totalement
oubli{\'e}e, {\`a} l'exception d'une br{\`e}ve remarque de Cohn dans
\cite{Cohn2} qui n'a pas attir{\'e} l'attention. La m{\'e}thode a
{\'e}t{\'e} r{\'e}invent{\'e}e par Kuhn en 1955 (\cite{Kuhn55}), sous
le nom de {\it m{\'e}thode hongroise}, en hommage au math{\'e}maticien
hongrois Egerv{\'a}ry ({\it cf.} \cite{Schrijver05} pour plus de
d{\'e}tails historiques).

Nous donnons ici une version de la borne de Jacobi dans le cas
quasi r{\'e}gulier, dans le cadre de la th{\'e}orie des diffi{\'e}t{\'e}s
({\it cf.} \cite{vinogradov,zharinov}).

\section{Diffi{\'e}t{\'e}s}

\begin{definition} Soit $I$ un ensemble d{\'e}nombrable, on munit $\R^{I}$
  de la topologie la plus grossi{\`e}re rendant pour tout $i_{0}\in I$ la
  projection $\pi_{i_{0}}:(x_{i})_{i\in I}\mapsto x_{i_{0}}$
  continue. On appelle {\em diffi{\'e}t{\'e}} un ouvert de $\R^{I}$ pour cette
  topologie, muni d'une d{\'e}rivation
$\delta=\sum_{i\in I} c_{i}(x){\partial\over\partial x_{i}}$,
o{\`u} les $c_{i}$ sont des applications $\cC^{\infty}$ ne d{\'e}pendant que d'un
{\em nombre fini} de coordonn{\'e}es. 
On note $\cO(V)$ l'anneau de telles applications sur la diffi{\'e}t{\'e} $V$.
\end{definition}

\begin{definition} On appelle morphisme de diffi{\'e}t{\'e}s une application
  $\phi:U_{\delta_{1}}\mapsto V_{\delta_{2}}$, d{\'e}finie par des
  fonctions $\cC^{\infty}$ ne d{\'e}pendant que d'un nombre fini de
  coordonn{\'e}es et telle que $\phi\circ\delta_{1}=\delta_{2}\circ\phi$.
\end{definition}
Ces d{\'e}finitions diff{\`e}rent de celles habituellement donn{\'e}es, o{\`u} une
diffi{\'e}t{\'e} est d{\'e}finie non par une d{\'e}rivation, mais par l'espace
vectoriel qu'elle engendre. Il n'en r{\'e}sulte pas de changement de fond,
mais quelques simplifications dans l'expos{\'e} de nos r{\'e}sultats.

\begin{definition} On appellera extension de diffi{\'e}t{\'e} et l'on notera
  $V/U$ un couple de diffi{\'e}t{\'e}s muni d'une projection
  $\pi:V\mapsto U$ surjective qui est un morphisme de
  diffi{\'e}t{\'e}s.  On appellera morphime d'extensions une
  application $\phi:V_{1}/U\mapsto V_{2}/U$ qui est un morphisme de
  diffi{\'e}t{\'e}s de $V_{1}$ dans $V_{2}$ tel que
  $\pi_{2}\circ\phi=\pi_{1}$.
\end{definition}

Cette d{\'e}finition {\it ad hoc} a pour but principal de correspondre
g{\'e}om{\'e}triquement aux extensions de corps diff{\'e}rentiels qui seraient
consid{\'e}r{\'e}es dans un cadre alg{\'e}brique.

\begin{definition}
Un ordre $\preceq$ sur les d{\'e}riv{\'e}es est dit admissible si
i) $u\prec u'$ et
ii) $u\preceq v$ implique $u'\preceq v'$.
\end{definition}

\begin{definition}\label{extension-triviale}
On appelle {\em extension triviale} de dimension diff{\'e}rentielle $m$ et l'on
note $T^{m}/U_{\delta_{0}}$ 
l'extension dont la d{\'e}rivation est donn{\'e}e dans des
coordonn{\'e}es $x_{i}^{(k)}$ par $\delta_{0}+\delta_{1}$ avec
$\delta_{1}=\sum_{i=1}^{m}\sum_{k=0}^{\infty} x_{i}^{(k+1)}\partial/\partial x_{i}^{(k)}$.
\end{definition}

\begin{definition}\label{forme-normale}
On dit qu'un syst{\`e}me $g$ d{\'e}finit une extension $V/U$ comme
sous-extension de l'extension triviale $T^{m}/U$ au voisinage de
$(v,u)$ s'il existe un ouvert $\cO\ni(v,u)$ de $T^{m}/U$ tel que $\cO\cap
V/U=\{(x,y)\in\cO|\forall k\in\N\forall i\in [1,r]
g_{i}^{(k)}(x,y)=0\}$. 

Une extension est dite {\em de type fini} si tout point admet un
voisinage isomorphe {\`a} une sous-extension de l'extension triviale.

Un syst{\`e}me de la forme $x_{i}^{(\gamma_{i})}=f_{i}(x)$, o{\`u} $f_{i}$ ne
d{\'e}pend que de d{\'e}riv{\'e}es des $x_{i}$ inf{\'e}rieures {\`a}
$x_{i}^{(\gamma_{i})}$ pour un ordre admissible, est appel{\'e} une {\em
forme normale} de l'extension qu'il d{\'e}finit dans les coordonn{\'e}es $x$.
\end{definition}

\begin{definition}\label{ordre-dimension}
Soit $V/U$ une extension de diffi{\'e}t{\'e}s,
d{\'e}finie comme sous extension de $T^{m}/U$ par une forme normale
$x_{i}^{(\gamma_{i})}=f_{i}(x)$, $1\le i\le r$.
On appelle {\em dimension differentielle} de l'extension
$V/U$ le nombre $m-r$ et {\em ordre de l'extension dans les coordonn{\'e}e
$x$ au voisinage de $(v,u)$} l'ordre maximal, g{\'e}n{\'e}riquement atteint, d'une
sous diffi{\'e}t{\'e} de dimension $0$, d{\'e}finie par des
{\'e}quations d'ordre nul en les $x_{i}$, au voisinage de $(v,u)$. 
\end{definition}
Si l'extension est de dimension $0$, son ordre ne d{\'e}pend pas des
coordonn{\'e}es et vaut $\sum_{i=1}^{n} \gamma_{i}$, mais sinon elle
peut d{\'e}pendre des coordonn{\'e}es et du point choisi.

\section{L'algorithme de Jacobi}

Nous allons d{\'e}crire tr{\`e}s bri{\`e}vement l'algorithme de
Jacobi, permettant de d{\'e}terminer, en notant $S_{r,m}$ l'ensemble des
injections de $[1,r]$ dans $[1,m]$, le maximum $\max_{\sigma\in S_{r,m}}
\sum_{i=1}^{r} a_{i,\sigma(i)}$.

\begin{definition} Soit $(a_{i,j})$ une matrice $m\times m$, on
  appelle {\em maximum} un terme sup{\'e}rieur ou {\'e}gal {\`a} tous
  les {\'e}l{\'e}ments de la m{\^e}me colonne, {\em maxima
  transversaux} des maxima situ{\'e}s dans des lignes et des colonnes
  toutes diff{\'e}rentes. On appelle {\em canon}
  un $m$-uplet d'entiers $(\lambda_{j})$ tels que la matrice
  $(a_{i,j} + \lambda_{j})$ poss{\`e}de $m$ maxima transversaux.
\end{definition}

\begin{lemme} Soit $(\ell_{j})$ et $(\lambda_{j})$ deux canons, alors
     $(\min(\ell_{j},\lambda_{j}))$ est un canon.
\end{lemme}
Soient $I_{1}$ (resp.\ $I_{2}$) l'ensemble des indices tels que
$\ell_{j}\le\lambda_{j}$ (resp.\ $\ell_{j}>\lambda_{j}$) et
$\sigma_{\alpha}\in S_{m,m}$, $\alpha=1,2$ telles que les
{\'e}l{\'e}ments d'indices $i,(\sigma_{\alpha}(i))$ constituent des
syst{\`e}mes de maxima transversaux pour les deux canons. Si
$j=\sigma_{1}(i_{1}\in I_{1})=\sigma_{2}(i_{2}\in I_{2})$, les
{\'e}l{\'e}ments d'indices $(i_{1},j)$ et $(i_{2},j)$ des deux canons
sont maximaux, et donc {\'e}gaux ce qui implique
$\ell_{i_{2}}=\lambda_{i_{2}}$ et donc $i_{2}\in I_{1}$, une
contradiction. Les images de $I_{1}$ et $I_{2}$ respectivement par
$\sigma_{1}$ et $\sigma_{2}$ sont donc disjointes et la r{\'e}union
des {\'e}l{\'e}ments d'indices $(i,\sigma_{\alpha}(i))$ pour $i\in
I_{\alpha}$ forme donc un syst{\`e}me de maxima transversaux.

L'algorithme d{\'e}crit par Jacobi dans \cite{Jacobi1} permet de
calculer en temps polynomial le canon minimal d'une matrice carr{\'e}e. Si $r=m$, on en
d{\'e}duit directement le maximum $J_{V,(v,u)}$. Sinon, on se ram{\`e}ne au cas d'une
matrice carr{\'e}e en compl{\'e}tant la matrice des ordres par des
lignes de $0$ (qui s'interpr{\`e}tent comme les ordre des
{\'e}quations g{\'e}n{\'e}riques de la d{\'e}finition
\ref{ordre-dimension}). On trouvera dans \cite{Munkres57} une version de la
m{\'e}thode hongroise et une {\'e}tude de sa complexit{\'e}.

\section{La borne}

\begin{definition}
On dit qu'un syst{\`e}me $g_{i}$ est {\em quasi r{\'e}gulier} en un point
de $V$ si pour tout $s\in\N$ la matrice jacobienne $J_{s}(g)$ des
applications d{\'e}riv{\'e}es $g_{i}^{(\ell)}$, $0\le \ell\le s$ par
rapport {\`a} toutes les d{\'e}riv{\'e}es des $x_{j}^{(k)}$ dont elles
d{\'e}pendent est en ce point de rang maximal {\'e}gal {\`a} $r(s+1)$.
\end{definition}

\begin{definition}
Soit $g(x,y)=0$ un syst{\`e}me de $r$ equations en $m$ inconnues
d{\'e}finissant $V$ au voisinage d'un point $(v,u)$ de $T^{m}/U$. Par
convention, on notera $\ord_{V,(v,u),x_{j}} g_{i}$ le plus grand
entier $k$ tel que $\partial g_{i}/\partial x_{j}^{(k)}\neq0$ sur tout
voisinage de $(v,u)$ dans $V$, ou $-\infty$ si $\partial
g_{i}/\partial x_{j}^{(k)}=0$ pour tout $k\in\N$ sur un voisinage de
$(v,u)$ dans $V$\footnote{Cette convention correspond {\`a} la borne
de Jacobi stricte.  La preuve de Lando n'est valable que pour la borne
faible, qui consid{\`e}re que l'ordre est nul dans ce cas.}. Soit
$a_{i,j}:=\ord_{V,(v,u),x_{j}} g_{i}$. Compl{\`e}tons
{\'e}ventuellement cette matrice de lignes de $0$ afin de la rendre
carr{\'e}e si $r<m$. Soit $\lambda_{i}$ un canon associ{\'e}.

 On posera
$\Lambda=\max_{i} \lambda_{i}$, $\alpha_{i}=\Lambda -\lambda_{i}$ et
$\beta_{j}=\max_{i} a_{i,j}-\alpha_{i}$.  La {\em matrice jacobienne
tronqu{\'e}e} $\nabla_{V,(v,u)}$ du syst{\`e}me $g$ est la matrice
$(\partial g_{i}/\partial x_{j}^{(\alpha_{i}+\beta_{j})})$.  On
appelle {\em ordre de Jacobi} un ordre tel que
$x_{j_{1}}^{(k_{1})}<x_{j_{2}}^{(k_{2})}$ si
$k_{1}-\beta_{j_{1}}<k_{2}-\beta_{j_{2}}$.  Les {\em d{\'e}riv{\'e}es
principales} de $g_{i}$ sont les d{\'e}riv{\'e}es
$x_{j}^{(\alpha_{i}+\beta_{j})}$.
\end{definition}

\begin{theoreme} Soit $V/U$ une extension de type fini de dimension
  differentielle $0$. Supposons la d{\'e}finie comme sous extension
  d'un ouvert de $T^{m}/U$ par un syst{\`e}me de $n$ {\'e}quations
  $g_{i}(x)=0$. On appelle nombre de Jacobi strict de $V$ en $(v,u)$
  l'entier $J_{V,(v,u)}:= \max_{\sigma\in S_{r,m}}
  \sum_{i=1}^{r} a_{i,\sigma(i)}$.

(i) Si le syst{\`e}me $g$ est quasi r{\'e}gulier en un point $(v,u)$ de
  $V/U$, alors il existe un ouvert $\cO$ dont l'adh{\'e}rence est un
  voisinage de $(v,u)$ et dont tout point admet un voisinage
poss{\'e}dant une forme normale et d'ordre au plus $J_{V,(v,u)}$ dans
les coordonn{\'e}es $x_{i}$.

(ii) Si $\nabla_{V,(v,u)}$ est de rang maximal en $(v,u)$, alors $g$ est
quasi r{\'e}gulier au voisinage de $(v,u)$, $V/U$ y est de dimension
$m-r$, admet une forme normale et est d'ordre
$J_{V,(v,u)}$ dans les coordonn{\'e}es $x_{i}$.

(iii) Si le syst{\`e}me $g$ est quasi r{\'e}gulier en un point $(v,u)$
de $V/U$ et si le rang de $\nabla_{V,(v,u)}$ n'est pas maximal sur un
voisinage ouvert de $(v,u)$ dans $V$, il existe un ouvert $\cO$ dont
l'adh{\'e}rence est un voisinage de $(v,u)$ et dont tout point admet
un voisinage poss{\'e}dant une forme normale et d'ordre strictement
inf{\'e}rieur {\`a} $J_{V,(v,u)}$ dans les coordonn{\'e}es $x_{i}$.
\end{theoreme}

\textsc{Preuve.} ---
Celle-ci suit le sch{\'e}ma d{\'e}crit par Jacobi dans \cite{Jacobi2} pour le
calcul d'une forme normale. On remarque d'abord que
$\ord_{V,x_{j}}g_{i}^{(s)}=\ord_{V,x_{j}}g_{i}+s$, en utilisant la
formule $[\partial/\partial
x_{j}^{(k+1)},\delta]=\partial/\partial x_{j}^{(k)}$, o{\`u} $\delta$ est
la d{\'e}rivation dans $T^{m}$.

(ii) Supposons que $\nabla_{V,(v,u)}$ soit de rang maximal en
$(v,u)$. On peut r{\'e}ordonner les
{\'e}quations $g_{i}$ de sorte que la suite $\alpha_{i}$ soit
croissante, et les $x_{i}$ afin que les mineurs principaux
constitu{\'e}s des $i$ premi{\`e}res lignes et colonnes de $\nabla_{V,(v,u)}$
aient un d{\'e}terminant $D_{i}$ non nul pour tout $1\le i\le r$. La
matrice jacobienne du syst{\`e}me $G_{s}:=\{g_{i}^{(k-\alpha_{i})}|0\le
k\le s\}$ contient une sous matrice carr{\'e}e de taille
maximale obtenue en d{\'e}rivant chaque {\'e}quation $g_{i}^{(k-\alpha_{i})}$ par
rapport aux d{\'e}riv{\'e}es principales 
$x_{j}^{(\beta_{j}+k)}$ telles que
$\alpha_{j}\le k$. Celle-ci est triangulaire par blocs et
son d{\'e}terminant est {\'e}gal {\`a} $\prod_{k=0}^{s} D_{\max\{i|
\alpha_{i}\le k\}}$. Cette matrice est donc pour tout $s$ de rang
maximal, ce qui entra{\^\i}ne la quasi-r{\'e}gularit{\'e} de $g$ au point
$(v,u)$.

En appliquant le th{\'e}or{\`e}me des fonctions implicites, on
obtiendra alors sur un voisinage de $(v,u)$ une relation de la forme
$x_{i}^{(\alpha_{i}+\beta_{i})}=f_{i}(v)$, $1\le i\le r$ o{\`u}
$f_{i}$ ne d{\'e}pend pas de $x_{j}^{(\alpha_{i}+\beta_{i})}$ pour
$j\le i$ et $\partial f_{i}/\partial x_{j}^{(k)}$ est nul sur $V$ pour
$k\ge \alpha_{j}+\beta_{j}$. Comme $g$ d{\'e}finit $V$, pour $s$ assez
grand, la fonction $f_{i}$ obtenue {\`a} partir du syst{\`e}me $G_{s}$
ne d{\'e}pendra que de d{\'e}riv{\'e}es inf{\'e}rieures {\`a}
$x_{i}^{(\alpha_{i}+\mu_{i})}$ pour l'ordre de Jacobi tel que
$x_{j_{1}}^{s+\beta_{j_{1}}}>x_{j_{2}}^{s+\beta_{j_{2}}}$ si
$j_{1}<j_{2}$. L'ordre dans les coordonn{\'e}es $x_{i}$ est alors
$\sum_{i=1}^{m} \alpha_{i}+\beta_{i}=J$ si $r=m$ et la
diffi{\'e}t{\'e} est de dimension $0$. En dimension positive, on se
ram{\`e}ne {\`a} ce cas en ajoutant $m-r$ {\'e}quations
g{\'e}n{\'e}riques d'ordre $0$, c'est {\`a} dire telles que la matrice
jacobienne tronqu{\'e}e demeure de rang maximal.

Pour (iii) et (i), on proc{\`e}de par r{\'e}currence sur $J_{V,(v,u)}$. Si
$J=0$, la quasi-r{\'e}gularit{\'e} implique que $\nabla_{V,(v,u)}$ est
de rang maximal et on utilise (ii) pour montrer (i); (iii) est alors
trivial puisque son hypoth{\`e}se est impossible.

(iii) Supposons que le syst{\`e}me soit quasi r{\'e}gulier et que
$\nabla_{V,(v,u)}$ s'annule identiquement sur un voisinage de
$(v,u)$. Supposant les {\'e}quations ordonn{\'e}es par $\alpha_{i}$
croissant, soit $k$ le plus petit indice tel que les $k$
premi{\`e}res lignes de $\nabla_{V,(v,u)}$ soient d{\'e}pendantes. On
peut proc{\'e}der pour le syst{\`e}me $g_{1}, \ldots, g_{i-1}$ comme
dans la preuve de (ii) et calculer au voisinage de $(v,u)$ une forme
normale de la diffi{\'e}t{\'e} qu'il d{\'e}finit. Les {\'e}quations
$g_{i}^{\ast}:=x_{i}^{(\alpha_{i}+\beta_{i})}-f_{i}(v)=0$ sont
{\'e}quivalentes {\`a} $g_{i}$ pour $i<k$. On peut alors remplacer
$g_{k}$ par une {\'e}quation {\'e}quivalente en substituant aux
d{\'e}riv{\'e}es principales $x_{j}^{(\alpha_{k}+\beta_{j})}$, $1\le
j<k$ les expressions $f_{j}^{(\lambda_{j}-\lambda_{k})}$. On obtient
ainsi une nouvelle {\'e}quation $g_{k}^{\ast}$ dont les
d{\'e}riv{\'e}es $x_{i}^{(\alpha_{i}+\beta_{i})}$, $i\le k$ ont {\'e}t{\'e}
{\'e}limin{\'e}es et telle que le syst{\`e}me $g_{1}^{\ast}, \ldots,
g_{k}^{\ast}, \ldots, g_{r}$ soit quasi r{\'e}gulier, {\'e}quivalent {\`a}
$g_{k}$ et tel que $\partial g_{k}/\partial
x_{i}^{(\alpha_{i}+\beta_{i})}$ s'annule sur $V$ au voisinage de
$(v,u)$. Le nouveau syst{\`e}me a donc un nombre de Jacobi strict
$J^{\ast}$ strictement inf{\'e}rieur {\`a} $J_{V,(v,u)}$, et on utilise (i) par
r{\'e}currence.

(i) Supposons que $|\nabla_{V,(v,u)}|=0$ au voisinage de $(v,u)$, en
utilisant (iii) on se ram{\`e}ne {\`a} un syst{\`e}me {\'e}quivalent
avec une borne de Jacobi strictement inf{\'e}rieure. Sinon, soit
$\cO'$ l'ouvert o{\`u} $|\nabla_{V,(v,u)}|\neq0$, $\cO''$
l'int{\'e}rieur du ferm{\'e} o{\`u} $|\nabla_{V,(v,u)}|=0$.
En utilisant (ii) sur $\cO'$ et (iii) sur $\cO''$, on
obtient (i) sur $\cO'\cup\cO''$ et $(v,u)$ appartient {\`a}
l'adh{\'e}rence de cet ouvert.

\section{Conclusion}

Dans le cadre des diffi{\'e}t{\'e}s, il est naturel de se restreindre {\`a} des
syst{\`e}mes quasi r{\'e}guliers. La borne de Jacobi fournit alors un r{\'e}sulat
th{\'e}orique meilleur que la \og borne de B{\'e}zout{\fg} diff{\'e}rentielle
$\sum_{j=1}^{m} \max_{i=1}^{r} a_{i,j}$, ce qui permet d'affiner les
r{\'e}sultats o{\`u} celle-ci intervient.

Les syst{\`e}mes quasi r{\'e}guliers incluent le cas de composantes
consid{\'e}r{\'e}es par Ritt comme singuli{\`e}res
(\cite{Ritt50}). {\`A} ce titre, la preuve de Kondratieva va au dela
de la situation initialement d{\'e}crite par Jacobi, puisqu'elle
utilise un nombre de Jacobi d{\'e}fini {\`a} partir de
$\ord_{V,x_{j}}g_{i}$ qui peut {\^e}tre inf{\'e}rieur {\`a}
$\ord_{x_{j}}g_{i}$.


\end{document}